\documentclass[10pt,twoside]{article}
\usepackage{graphicx}
\usepackage{amsmath}
\usepackage{Latex-document}
\newenvironment{extrait}
{\begin{list}{}{\leftmargin 5mm \rightmargin0pt}%
\item\itshape}{\end{list}}

\newcommand{\beqn} {\begin{eqnarray*} }
\newcommand{\eeqn} {\end{eqnarray*} }

\newcommand{\beqnn} {\begin{eqnarray} }
\newcommand{\eeqnn} {\end{eqnarray} }

\newcommand{\bc} {\begin{center} }
\newcommand{\ec} {\end{center} }

\newcommand{\bit} {\begin{itemize} }
\newcommand{\eit} {\end{itemize} }

\newcommand{\bnum}{\begin{enumerate}}
\newcommand{\enum}{\end{enumerate}}

\def \hs{\hspace{.5cm}}

\def \lim{{\rm lim}}
 \def \Log{{\rm Log}}

\def \E{{\mathbf E}}
\def \Qr{{\mathbf Q}}

\def \demi{\frac{1}{2}}

\title{\bf Measuring and Hedging  \vskip -2mm
Financial Risks in Dynamical World  \vskip 6mm}
\author{Nicole El Karoui\vspace*{-0.5cm}\thanks{Centre de Math\'ematiques Appliqu\'ees, Ecole Polytechnique, 91 128
Palaiseau Cedex, France. E-mail: elkaroui@cmapx.polytechnique.fr}}
\date{\vspace{-8mm}}

\begin{document}
\maketitle \thispagestyle{first} \setcounter{page}{773}
\begin{abstract}\vskip 3mm
Financial markets have developed a lot of strategies to control
risks induced by market fluctuations. Mathematics has emerged as
the leading discipline to address fundamental questions in finance
as asset pricing model and hedging strategies. History began with
the paradigm of zero-risk introduced by Black \& Scholes stating
that any random amount to be paid in the future may be replicated
by a dynamical portfolio. In practice, the lack of information
leads to ill-posed problems when model calibrating. The real world
is more complex and new pricing and hedging methodologies have
been necessary. This challenging question has generated a deep and
intensive academic research in the 20 last years, based on
super-replication (perfect or with respect to confidence level) and optimization. In the interplay between theory and
practice, Monte Carlo methods have been revisited, new risk measures have been back-tested.
 These typical examples give some insights
 on how may be used mathematics in financial risk
management.


\end{abstract}

\vskip 12mm

\section{Introduction}\label{section 1}\setzero
\vskip-5mm \hspace{5mm}

Financial markets have become an important component
of people's life. All sorts of media
now provide us with a daily coverage on financial news from all
markets around the world.
At the same time, not only large institutions but also more and more
small private investors
are taking an active part in financial trading. In particular, the
e-business has led to an unprecedented increase
in small investors direct trading.
Given the magnitude of the potential impact a financial crisis can
have on the real side of the economy,
large, and even small, breakdowns have received particularly focal attention.
Needless to say, the rapid expansion
of financial markets  calls
upon products and systems
designed
to help investors to manage their financial risks. The financial risk business
now represents more than \$15 trillion annually in notional.
A large spectrum of simple contracts (futures, options,
swaps, etc.) or more exotic financial products
(credit derivatives, catastrophe bonds, exotic options, etc.)
are offered to private investors
who use them to transfer financial risks to
specialized financial institutions in exchange for suitable compensation.
A classic example is the call option, which provides a protection in case of a
large increase
in the underlying asset price
\footnote{A Call option provides its buyer with the
right
(and not the obligation) to purchase the risky asset at a
pre-specified price (the exercise price) at or before a pre-specified
date in the future (maturity date). The
potential gain at maturity
 can therefore be written as $(X_T-K)^+$, where
$X_T$ denotes the value of the underlying asset at time $T$}. More generally, a derivative contract is an asset
that delivers a payoff $H(\omega)$ at maturity date, depending upon the scenario $\omega$.

As argued by Merton \cite{Mert}, the development of the financial
risk industry would not have been possible without the support of
theoretical tools. Mathematics has emerged as the leading
discipline to address fundamental questions in financial risk
management, as asset pricing models and hedging strategies, based
on daily (infinitesimal) risk management. Mathematical finance,
which relates to the application of the theory of probability and
stochastic processes to finance, in particular the Brownian motion
and martingale theory, stochastic control and partial differential
equations, is now a field of research in its own right.
\section{The Black \& Scholes paradigm of zero-risk}
\label{section 2}\setzero
\vskip -5mm \hspace{5mm }

It is surprising that the starting point of financial industry
expansion is ``the Brownian motion theory and It\^o stochastic
calculus", first introduced in finance by Bachelier in this PhD
thesis (1900, Paris), then used by Black, Scholes and Merton in
1973. Based on these advanced tools, they develop the totally new
idea in the economic side that according to an optimal dynamic
trading strategy, it is possible for option seller to deliver the
contract at maturity without incurring any residual risk. At this
stage, any people not familiar with stochastic analysis (as the
majority of traders in the bank) may be discouraged. As Foellmer
said in  Bachelier Congress 2000, it is possible to reduce
technical difficulties,
and to develop arguments which are essentially probability-free.\\[1mm]
  {\sc A dynamical uncertain world }
  \\[1mm]\indent
 The uncertainty is modelled via a family
$\Omega$ of scenarios $\omega$, i.e. the possible trajectories of the asset prices in the future. Such paths are
described as positive continuous functions $\omega$ with coordinates $X_t=\omega(t)$, such that the continuous
quadratic variation exists: $[X]_t(\omega)=\lim_n\sum_{t_i\leq t,t_i\in D_n}(X_{t_{i+1}}-X_{t_i})^2$ along the
sequence of dyadics partitions $D_n$. The pathwise version of stochatic calculus yields to It\^o's formula, \bc $
f(t,X_t)(\omega)=f(0,x_0)+\int_0^tf_x'(s,X_s)(\omega)\,dX_s(\omega)+
\int_0^t\,f_t(s,X_s)(\omega)\,dt$ \\[2mm]
$+ \int_0^t \frac 12 f''_{xx}(s,X_s)(\omega) d[X]_s(\omega). $
 \ec
The second integral is well defined as a Lebesgue-Stieljes
integral, while the first exists as It\^o's integral, defined as
limit of non-anticipating Riemann sums, (where we put
$\delta_t=F'_x(t,X_t)$),
$\sum_{t_i\leq t,t_i\in D_n}\xi_{t_i}(\omega)(X_{t_{i+1}}-X_{t_i})$.\\[1mm]
\indent From a financial point of view, the It\^o's  integral may be interpreted as the cumulative {\em gain
process} of trading strategies: $\delta_t$ is the number of the shares held at time $t$ then the increment in the
Riemann sum is the price variation over the period. The non anticipating assumption corresponds to the financial
requirement that the investment decisions are based only on the past prices observations. The residual wealth of
the trader is invested only in cash, with yield rate (short rate) $r_t$ by time unit. The {\em self-financing
condition} is expressed as:
 \bc
\vskip -1mm $ \label{eq: self} dV_t=r_t(V_t-\delta_t.X_t)dt
+\delta_t.dX_t = r_tV_tdt +\delta_t.(dX_t-r_tX_t\,dt),\hs V_0=z. $
\ec
 {\sc Hedging derivatives: a solvable target
 problem}\\[1mm]\indent
Let us come back to the problem of the trader having to pay at
maturity $T$ the amount $h(X_T)(\omega)$ in the scenario $\omega$
($(X_T(\omega)-K)^+$ for a Call option). This target has to be
hedged (approached) in all scenarios by the wealth generated by a
self-financing portfolio. The "miraculous" message is that, in
Black \& Scholes  world,
 a perfect hedge is possible
and easily computable, under the additional assumption: the short rate is deterministic and the quadratic
variation is  absolutely continuous $d[X]_t=\sigma(t,X_t)\,X_t^2 dt$. The (regular) strictly positive function
$\sigma(t,x)$ is a key parameter in financial markets, called the {\em local
  volatility}.

Looking for the wealth as a function
$f(t,X_t)$, we see that, given It\^o's formula and self-financing condition,
\bc
$
df(t,X_t)=f'_t(t,X_t)dt+f'_x(t,X_t)dX_t+ \frac 12
f''_{xx}(t,X_t)X_t^2\sigma^2(t,X_t)\,dt $\\[2mm]
$=f(t,X_t)r_tdt+\delta(t,X_t)\big(dX_t-X_t r_t\,d
t\big)$
\ec
By identifying the $dX_t$ terms (tanks to assumption
$\sigma(t,x)>0$), $\delta(t,X_t)=f'_x(t,X_t)$, and
$f$ should
 be solution of the following partial differential equation,
Pricing PDE in short,
\bc
$f'_t(t,x)+\frac 12 f''_{xx}(t,x)x^2 \sigma^2(t,x) + f'_x(t,x)x
r_t-f(t,x)r_t =0, \> f(T,x)=h(x)\hs (2.1)$ \ec The derivative
price at time $t_0$ must be $f(t_0,x_0)$, if not, it is easy to
generated profit without bearing any risk (arbitrage). That is the
{\em rule of the unique price}, which holds in a liquid market.

The PDE's fundamental solution $q(t,x,T,y)$ ($h(x)=\delta_{y}(x)$) may
 be interpreted in terms of Arrow-Debreu ``states prices'' density,
 introduced in 1953 by these Nobel Prize winners
 for a purely theoretical economical point of view and by completely
 different arguments. The pricing rule becomes: $f(t,x)=\int h(y)q(t,x,T,y) dy$.
 $q$ is also called {\em pricing kernel}. When $\sigma(t,x)=\sigma_t$, the
  pricing kernel is the {\em log-normal density},
deduced from the Gaussian distribution by an explicit change of
variable. The closed formula for Call option price is the famous
\footnote{In the Black-Scholes model with constant
  coefficients, the Call option price
$C^{BS}(t,x,K,T)$ is given via the Gaussian cumulative function
${\cal N}(z)=\int_{\infty}^z\frac{1}{\sqrt{2\pi}}
e^{-\frac{y^2}{2}}dy$,and $\theta=T-t$,
\beqn\label{eq:BLS}
\left\{
\begin{array}{clllcr}
C^{BS}(t,x,t+\theta,K)=x\,{\cal N}\big(d_1(\theta,x/K)\big)-
K\,e^{-r\,\theta}{\cal N}\big(d_0(\theta,x/K)\big) \\[2mm]
d_0(\theta,x/K)=\frac{1}{\sigma \sqrt{\theta}}\Log
\left(\frac{x}{Ke^{-r\, \theta}}\right)-\frac{1}{2}\sigma
\sqrt{\theta},\hs
d_1(\theta,x/K)=d_0(\theta,x/K)+\sigma \sqrt{\theta}
\end{array}
\right.
\eeqn
 Moreover
 $\Delta(t,x)=\partial_xC^{BS}(t,x,t+\theta,K)={\cal N}\big(d_1(\theta,x/K)\big)$.}
Black and Scholes formula, which is known by any practitioner in
finance. The impact of this methodology was so
important that Black (who already died), Scholes and Merton received the Nobel prize for
economics in 1997.

In 1995, B.Dupire \cite{Dupi} give a clever formulation for the dual PDE (one
dimensional in state variable) satisfied by $q(t,x,T,y)$ in the
variables $(T, y)$. If $C(T,K)$ is the Call price
with parameters $(T, K)$ when market conditions are $(t_0,x_0)$, then
\bc
$ C'_T(T,K)=\frac 12 \sigma^2(T,K)\,K^2 C''_{KK}(T,K)
-rK\,C'_K(T,K),\>C(t_0,x_0)=(x_0-K)^+
$
\ec
In short, if $r_t=0$, $
\hs C'_T(T,K)=\frac 12 \sigma^2(T,K)\,K^2 C''_{KK}(T,K)$

\section{Model calibration and Inverse problem}
\label{section 3}\setzero
\vskip-5mm \hspace{5mm }

In pratice, the main problem is the model calibration. Before
discussing that, let me put the problem in a more classical
framework. Following P. L\'evy, asset price dynamics may be
represented through a Brownian motion via stochastic differential
equation (SDE)
\bc
$
dX_t=X_t(\mu(t,X_t)dt+\sigma(t,X_t)dW_t),\>\>X_{t_0}=x_0
$\ec
where the Brownian motion $W$ may be viewed as a standardized
Gaussian noise with independent increments.

The local expected return
$\mu(t,X_t)$ is a trend parameter appearing for the first time in our
propose. That is a {\em key point}
in financial risk management. Since this parameter does not appear in
the Pricing-PDE, the Call price does not depend on the market trend.
It could seem surprising, since the first motivation of this
financial product is to hedge the purchaser against underlying rises.
By using dynamical hedging strategy, the trader (seller) may be
also protected against this unfavorable evolution.
For a statistical point of view, this point is very important,
because this parameter is very difficult to
estimate.

In the B \& S model with constant parameters, the volatility square is
the variance by time unit of the log return
$\ln(X_t)-\ln(X_{t-h})=R_t$.
If the only available information
 is given by asset price historical data, the statistical
estimator to be used is the empirical variance, computed on a
more or less long time period,
$\widehat
\sigma^2=\frac{1}{N-1}\sum_{i=0}^{N-1}(R_{t_i}-\overline{R_N})^2$,
where
$\overline{R_N}=\frac{1}{N}\sum_{i=0}^{N-1}R_{t_i}$.
This estimator is called {\em
  historical volatility}.

However, traders are reserved in using this estimator. Indeed, they argue that financial markets are not
``statistically'' stationary and that past is not enough to explain the future. When it is possible, traders use
additional information given by quoted option prices and translate it into volatility parametrization. The {\em
implied volatility} is defined as:
$C^{obs}(T,K)=C^{BS}(t_0,x_0,T,K,\sigma^{imp})$.

Moreover, the $\Delta$ of the replicating portfolio is
$\Delta^{imp}=\partial_x C^{BS}(t_0,x_0,T,K,\sigma^{imp})$.

This strategy is used dynamically by defining the implied volatility and the
associated $\Delta$ at any renegotiation dates.
It was this specific use based on the hedging strategy that may
explain the Black \& Scholes formula success.

This attractive methodology has been appeared limited: observed implied volatilities are depending on the option
parameters (time to maturity and exercise price) (implied volatility surface) in complete contradiction with B\&S
model assumptions. In particular, the market quotes large movements (heavy tail distribution) higher than in the
log-normal framework. The first idea to take into account this empirical observation is to move to a model with
local volatility $\sigma(t,x)$. The idea is especially attractive since the Dupire formula (2.2) gives a simple
relation between ``quoted Call option prices'' and local volatility: $\sigma^2(K,T)=
2C'_T(T,K)(K^2C''_{KK}(T,K))^{-1}$. The local volatility is computable from a continuum of {\em coherent} (without
arbitrage) observed quoted prices. Unfortunately, the option market is typically limited to a relatively few
different exercize prices and maturities; a naive interpolation yields to irregularity and instability
of the local volatility.\\[1mm]
{\sc Ill-posed inverse problem}\\[1mm]\indent
The problem of determining local volatility  can be viewed as a
function approximation non linear problem from a finite data set.
The data set is the value $C_{i,j}$ of the solution at $(t_0,x_0)$
of Pricing-PDE with boundary conditions $h_{i,j}(x)=(x-K_{i,j})^+$
at maturity $T_i$. Setting the problem as PDE's inverse problem
yields to more robust calibration methods. These ideas appear for
the first time in finance in 1997 \cite{LaOs}, but the problem is
not classical because of the strongly non linearity between option
prices and local volatility;
the data set is related to a single given initial condition.

Prices adjustment is made through a least square minimization
 program, including a penalization term  related to the local volatility regularity.
\bc
$
G(\sigma)=
\sum_{i,j}\omega_{i,j}
\big(f(t_0,x_0,h_{i,j},T_i,\sigma(.,.))-C_{i,j}^{Obs}\big)^2,$\\[1mm]
$J(\alpha,\sigma)=\alpha ||\nabla \sigma||^2+ G(\sigma)\rightarrow \min_{\sigma}.$ \ec Existence and uniqueness of
solution is only partially solved \cite{BeBF}. Using  large deviation theory, the asymptotic in small time of
local volatility is expressed in terms of implied volatility: $ \sigma^{\rm
implied}(K,t_0)^{-1}=\ln(\frac{K}{x_0})^{-1} \int_{x_0}^K\frac{d\xi}{\xi\sigma(\xi,t_0)}
$.\\[1mm]\indent
Avellaneda \& alii \cite{AFHS} have used another penalization criterion
based on a stochastic control approach ; the control is the volatility
parameter itself constrained to be very close to a prior
volatility ($\eta(\sigma)
=|\sigma(t,x)-\sigma_0(t,x)|^2$ for instance). The gradient criterion
is replaced by $K(\sigma)=U(t_0,x_0,\sigma)$
 where $U(T,x,\sigma)=0$ and
\bc
$
U'_t(t,x)+\frac 12 \sigma^2(t,x)\,x^2 U''_{xx}(t,x)+
rx\,U'_x(t,x)-U(t,x)
+{\mathbf \eta(\sigma(t,x))}=0.
$
\ec
\section{Portfolio, duality and incomplete market}
\label{section 4}\setzero
\vskip-5mm \hspace{5mm }

In the previous framework, options market may be entirely explained by underlying prices. In economic theory, it
corresponds to {\em market efficiency}: a security price contains all the information on this particular security.
In option world, the observed statistical memory of historical volatility leads naturally to consider stochastic
volatility models with specific uncertainty
 \bc
$ dX_t=X_t\big(\mu(t,X_t,Y_t)dt + \sigma(t,X_t,Y_t) dW^1_t\big),\hs dY_t=
\eta(t,X_t,Y_t)+\gamma(t,X_t,Y_t)\,dW^{2}_t $ \ec where $dW^1$ and $dW^2$ are two correlated Brownian motions.
$\gamma$ is the volatility of the volatility. What does it change ? In fact, everything ! Perfect replication by a
portfolio is not possible any more ; the notion of unique price does not exist any longer... But, such a situation
is often the general case. What kind of answer may we bring to such a problem?\\
{\bf   Super-replication and Robust Hedging}\\[2mm]\indent
The option problem is still a target problem $C_T$, to be replicated
by a portfolio \\[1mm]
$V_T(\pi,\delta)=\pi+\int_0^T\sum_{i}\delta^i_s
dX_t^i$ depending on market assets $X^i$.
Constraints (size, sign...) may be imposed on investment decisions
$(\delta^i_t)$ \footnote{For the sake of simplicity, interest
  rates are assumed to be null}.
Let ${\cal V}_T$ be the set of all derivatives, replicable at time
$T$ by an admissible portfolio. Their price at $t_0$  is the value
of
 their replicating portfolio.
\begin{extrait}{\rm
Super-replicating $C_T$ is finding the smallest derivative
 $\widehat C_T\in {\cal V}_T$ which is greater than $C_T$ in all scenarios.
The super-replication price is the price of such a derivative. The
 $\widehat C_T$ replicating portfolio is the $C_T$ {\em robust
 hedging}.}
\end{extrait}
There are several ways to characterise the super-replicating
portfolio:

1) Dynamic programming on level sets: this is the most direct (but
least recent) approach. This method proposed by Soner \& Touzi
\cite{SoTo} has led the way for original works in geometry by
giving a stochastic representation of a class of
mean curvature type geometric equations.

2) Duality: this second approach is based on the ${\cal V}_T$ ``orthogonal space'', a set ${\cal Q}_T$ of
martingale measures to be characterised. The super-replication price is given by \bc \vskip -1mm $\widehat
C_0=\sup_{Q\in {\cal Q}_T} \E_{Q}[C_T]. $
 \ec
We develop this last point, which is at the origin of many works.\\[2mm]
{\bf Martingale measures}\\[2mm]\indent
 The idea of introducing a dual theory based on probability measures
is due to Bachelier (1900), and above all to
Harisson \& Pliska (1987). The actual and achieved form is due to
Delbaen \& Schachermayer \cite{DeSc94} and to the very active
international group in Theoretical Mathematical Finance.

 A {\em martingale measure} is a probability such that:
 $\forall V_T\in {\cal V}_T,\>\E_{Q}[V_T]=V_0$.
Using simple strategies (discrete times, randomly chosen), this
 property is equivalent to prices of fundamental assets
 $(X^i_t)$ are $\Qr$-(local) martingales: the best $X^i_{t+h}$-estimated
 (w.r. to $\Qr$) given the past at time $t$ is
 $X^i_{t}$ itself. The financial game is fair with respect to martingale measures.

When ${\cal V}_T$  contains all possible (path-dependent) derivatives, the market is
said to be complete,
and the set of martingale measures is reduced to a unique element
$\Qr$, often called risk-neutral probability. This is the
case in the previous framework. Dynamics become $dX_t=X_t\, \sigma(t,X_t) dW^Q_t $
where $W^Q$ is a $\Qr$-Brownian motion.
This formalism is really efficient as it leads to the following
path dependent derivative {\em pricing rule}: $\widehat C_{0}=\E_{Q}(C_{T}).$

Computing the replicating
portfolio is more complex. In the diffusion case,
the price is a deterministic function of risk factors and the replicating portfolio
only depends on partial derivatives.
The general case will be mentioned
in the paragraph dedicated to Monte-Carlo methods.\\[1mm]
{\sc Incomplete market}\\[1mm]\indent
The characterization of the set ${\cal Q}_T$ is all the more delicate
so since
there are many different situations which may lead to market
imperfections (non-tradable risks, transaction costs ...).

An abstract
theory of super-replication (and more generally of portfolio
optimization under constraints) based on duality has been intensively developed.
The super-replicating price process is showed \cite{NEKQ},\cite{Kram}
to be a super-martingale with respect to any admissible martingale measure.
Hence, by the generalization of the Doob-Meyer representation, the
super-replicating
portfolio is the ``${\cal Q}_T$- martingale'' part
of the super-price Kramkov-decomposition.

Super-replication prices are often too expensive to be used in practice. However, they give an upper bound to the
set of possible prices. In the previously described stochastic volatility model, the super-replication price
essentially depends on possible values of stochastic volatility:
\begin{enumerate}
\item  If the set is $R^+$, then the super-replicating
derivative of $h(X_T)$ is $\widehat h(X_T)$
 where $\widehat h$
is the concave envelop of $h$ ; the replicating strategy is the
trivial one: buying $\widehat h'_x(x_0)$ stocks and holding them
till maturity. \vskip -0.5mm
\item If the volatility is bounded (up and down relatively to 0), the
super-replication price is a (not depending on $y$) solution of \bc $ {\widehat h}'_t(t,x)+\demi
\sup_y(\sigma^2(t,x,y){\widehat h}''_{xx}(t,x))=0,\hs \widehat h(T,x)=h(x).$ \ec When $h$ is convex, $\widehat
h(t,x)$ is convex and the super-replication price is the one calculated with the {\em upper volatility} (in $y$).
\end{enumerate}
\indent Calibration constraints may be easily taken into account
without modifying this framework. We only have to assume that the
terminal net cash flows of calibrating derivatives belong to
${\cal V}_T$ or equivalently we have to add linear constraints to
the dual problem:  $ \hs \widehat C^{{\rm cal}}_0=
\sup\{\E_{Q}(C_T)\>; Q\in {\cal Q}_T\,,
\E_{Q}((X_{T_i}-K_{ij})^+)= C_{ij}\}. $
\\[2mm]
{\bf   Risk measures}\\[2mm]\indent
When super-replicating is too expensive, the trader has
to measure his market risk exposure. The traditional
measure is the variance of the replicating error. But a new criterion, taking into
account extreme events,
is now used, transforming the risk management at both quantitative
and qualitative levels.\\[1mm]
{\sc Value at Risk}\\[1mm]\indent
The VaR criterion, corresponding to {\em the maximal level of losses}
acceptable with probability 95\%,
has taken a considerable importance for several years. Regulation
Authorities
have required a daily VaR computation of the
global risky portfolio from financial institutions.
Such a measure is important on the operational point of view,
as it affects the provisions a bank has to hold to face market risks.
VaR estimation (quantile estimation) and its links with extreme value theory
\cite{Embr} are widely debated in the market, just as by academics.

Moreover, a huge debate has been introduced by academics
\cite{ADEH} on the VaR efficiency as risk measure.
For instance, its non-additive property enables banks to play with
subsidiary creations. This debate has received an important echo
from the professional world,
 which is possibly planning to review this risk measure criterion. Sub-additive and coherent risk measures are an average estimation of
losses with respect to a probability family:
$\rho(X)=\sup_{Q\in {\cal Q}_T}\E_{Q}(-X)$.

This characterisation has recently been extended to convex risk
measures, by adding a penalization term depending on probability density
(entropy for instance).\\[1mm]
{\sc Risk measures and reserve price}\\[1mm]\indent
A trader willing to relax the super-replication assumption
is naturally thinking in terms of potential losses with given
probability (level confidence)
of his replicating error. It corresponds to the quantile hedging
strategy.
Other risk measures (quadratic, convex, entropic) may be used.
Optimization theory is coming back when looking for the smallest
portfolio, generating an acceptable loss.
The initial value of this portfolio is called the {\em reserve price}.
Mean-variance and entropic problems have now a complete solution \cite{Frit}.
 More surprisingly (because of non-convexity), this
also holds for the quantile hedging problem \cite{FoLe}.
All these results are in fact sub-products of portfolio optimization
in incomplete markets \cite{CvKa} or \cite{Duff}.

\section{New research fields}
\label{section 5}\setzero
\vskip-5mm \hspace{5mm }

\noindent{\bf   Monte-Carlo methods}\\[2mm]\indent
Dual version of super-replication problems,
just as new risk measures, underline the interest to compute very well
and quickly quantities such as $\E_{Q}(X)$ and more generally
$\sup_{Q\in {\cal Q}_T}\E_{Q}(X)$.
For small dimensional diffusions, these quantities may be
computed as the solution
of some linear PDE (for the expected value)
and non-linear PDE (for the sup).
However, the computational efficiency falls rapidly with the
dimension. That increases the interest for the so-called {\em probabilistic methods}.

The fundamental idea of Monte-Carlo methods is the computation of
$\E_{Q}(X)$ by simulation, i.e.
by drawing a large number
($N\simeq 10^5$)
of independent scenarios
$\omega^i$
 and taking the average value of the results
$\frac 1N\sum_{i=1}^N X(\omega^i).$
Of course, this method does not work very well when being too naive,
 but convergence may be accelerated by different techniques.

In the finance area, the important quantities are both the price and the
 sensitivities to different model parameters.
 Based on integration by parts,
 efficient methods have been developed to compute in a coherent manner
 prices and their
 derivatives \cite{LilaF2}.
  In the case of path-dependent options,
  the derivative is taken with respect
 to a Brownian motion perturbation (Malliavin Calculus).

Very original is the actual research
 on solving, by Monte-Carlo methods, optimization problem
 expressed as "sup" of a family of expected values
 (super-replication prices).
 These solutions are based on dynamic programming,
 enabling to turn a maximization in expected value into
 a pathwise maximization. The formulation in terms of Backward SDE's,
 introduced by Peng and Pardoux in 1987 and in
 finance\footnote{A BSDE solution is a couple of adapted processes
 $(Y_t,Z_t)$ such that $$-dY_t=f(t,Y_t,Z_t)dt-Z_tdW_t,\hs Y_T=C_T.$$
 Thanks to comparison theorem, ${\widehat Y}_t=\sup_{i\in I}Y^i_t$
 is the solution of the BSDE with driver
${\widehat f}(t,y,z)=\sup_{i\in I}f^i(t,y,z)$}
 in \cite{EPQPbsde},\cite{MaYo} well describes this effect.

In all cases, the problem is to compute conditional expectation
by Monte Carlo, using  the function
approximation theory, or more generally
random variable approximation in
in the Wiener space (chaos decomposition).\\[2mm]
{\bf  Problems related to the dimension, and statistical modelling}\\[2mm]\indent
Financial problems are usually multidimensional, but only few
 {\em liquid} financial products are depending on  multi-assets.
Even if the different market actors consider
 they can have a good knowledge of
 each individual asset behavior, the question
is now to find a multidimensional distribution
given each component distribution. This problem is a statistical one,
 known as the {\em copula theory}:
 Copula is a distribution function on $[0,1]^n$
 with identity function as  marginal. They are useful to give
 bounds to asset prices.
 Dynamically, the problem still to be solved is to find  the local volatility matrix of
 multidimensional diffusion given the "Dupire" dynamics of each coordinate.

High dimensional problems arise when computing bank portfolio VaR (with
 a number of observations less than that of risk factors), or
hedging derivatives depending on a large number of
 underlying assets. Main risk factors may be very different in a
 Gaussian framework or heavy tail framework (L\'evy processes)\cite{BoSoWA}.
Random matrix theory or other asymptotic tools may
 bring some new ideas to this question.

 By presenting the most important tools
 of the financial risk industry, I have voluntary left apart anything
 on financial asset statistical modelling,
 which may be the subject of a whole paper on its own.
 It is clear that the VaR criterion, the market imperfections are
 highly dependent
 on an accurate analysis of the real and historical world \cite{BaNi}.
 Intense and very innovating research is now developed (High-frequency data,
ARCH and GARCH processes, L\'evy processes with
 long memory, random cascades).

\section{Conclusion}
\label{section 6}\setzero \vskip-5mm \hspace{5mm }

As a conclusion, applied mathematicians have been highly questioned by
 problems coming from the financial risk industry.
 This is a very active world, rapidly evolving, in which theoretical
 thoughts
 have often immediate and practical fallout. On the other hand,
 practical
 constraints raise new theoretical problems.
 This paper is far from being an exhaustive view of the financial
 problems.
 It is more a subjective view conditioned by my own experience.
 Many exciting problems, from both theoretical and practical points
 of view,
 have not been presented. May active researchers in these fields
 forgive me.

\label{lastpage}

\end{document}